\documentclass[12pt,a4paper]{article}
\usepackage{amssymb,amsmath}
\usepackage{amsfonts}
\usepackage{a4}
\usepackage{tikz}
\usepackage{verbatim}
\usepackage{caption}
\usepackage{float}
\newcommand{\para}{\par\vspace{.25cm}}

\newtheorem{prop}{Proposition}
\newtheorem{theorem}{Theorem}
\newtheorem{lemma}{Lemma}
\newtheorem{cor}{Corollary}

\begin{document}
\baselineskip 18pt \title{\bf \ A generalization of strongly monomial groups}
\author{ Gurmeet K. Bakshi and Gurleen Kaur{\footnote {Research supported by NBHM(Ref No: 2/39(2)/2015/NBHM/R\&D-II/7424), Department of Atomic Energy, Government of India, is gratefully acknowledged.} \footnote{Corresponding author}} \\ {\em \small Centre for Advanced Study in
Mathematics,}\\
{\em \small Panjab University, Chandigarh 160014, India}\\{\em
\small email: gkbakshi@pu.ac.in and gurleenkaur992gk@gmail.com  } }
\date{}
{\maketitle}
\begin{abstract}  Olivieri, del R{\'{\i}}o and Sim{\'o}n defined strongly monomial groups and a significant result proved by them is the explicit description of the simple components of the rational group algebra $\mathbb{Q}G$ of a strongly monomial group $G$. In this paper, generalized strongly monomial groups are defined, which is a generalization of strongly monomial groups. Beside strongly monomial groups, an extensive list of important families of groups which are generalized strongly monomial is produced. A description of the simple components of $\mathbb{Q}G$, when $G$ is a generalized strongly monomial group, is also provided. Furthermore, the work by Jespers, Olteanu, del R{\'{\i}}o and Van Gelder on the construction of a subgroup of finite index in the group of central units of integral group ring of a strongly monomial group has also been generalized.\end{abstract}
{\bf Keywords} : rational group algebra, simple component, strongly monomial group, generalized strongly monomial group, integral group ring, central unit. \vspace{.25cm} \\
{\bf MSC2000 :} 16S34, 16K20, 16S35, 20C05, 16U60\section{Introduction}Throughout this paper, $G$ is a finite group and  $\mathbb{Q}G$ is the rational group algebra of $G$. One of the fundamental problem in group rings is to understand the primitive central idempotents and the simple components of $ \mathbb{Q}G$. The classical approach using character theory is computationally hard and does not yield any insight into the  structure of the simple components of $\mathbb{Q}G$. A new direction to this investigation was provided by Olivieri, del R{\'{\i}}o and Sim{\'o}n in \cite{OdRS04} where they introduced the notion of strong Shoda pairs and strongly monomial groups. A strong Shoda pair of $G$ is a pair of subgroups of $G$ satisfying some constraints. These pairs of subgroups form the main ingredient to define a strongly monomial group. A significant result proved in \cite{OdRS04} is the description of the simple components of $ \mathbb{Q}G$ when $G$ is strongly monomial. This research shifted the focus of group ring theorists to understand strong Shoda pairs and strongly monomial groups and to seek applications in understanding the unit group of integral group rings (see \cite{Ba},\cite{BK},\cite{BM},\cite{BM1},\cite{BMP},\cite{JOdRVG1},\cite{JOdRVG},\cite{JP},\cite{M}). In this paper, we generalize the notion of strong Shoda pairs and strongly monomial groups  and consequently generalize the main work of \cite{JOdRVG} and \cite{OdRS04}. \par
When $(H,K)$ is a strong Shoda pair of $G$, then any linear character $ \lambda $ on $H$ with kernel $K$ when induced to $G$ is irreducible and, in addition, $\lambda$ satisfy some constraints. In section 2, we generalize this concept and give the notion of a strong inductive source of $G$. By a strong inductive source of $G$, we mean an irreducible character of a subgroup of $G$ which induced to $G$ is irreducible, and it has the same constraints which the linear character $\lambda$ had in the definition of a strong Shoda pair of $G$. This perception allows us to define generalized strong Shoda pairs of $G$ and accordingly the generalized strongly monomial groups. Our motivation behind this generalization has been the properties of Shoda pairs constructed in \cite{BK}. We can immediately conclude from the definition  that all strongly monomial groups and subnormally monomial groups are generalized strongly monomial. Beside these groups, we show that the class of generalized strongly monomial groups contain the groups in $\mathcal{C}$. The class $\mathcal{C}$ consists of all finite groups $G$ such that all subgroups and quotient groups of subgroups of $G$ satisfy the following property: either it is abelian or it contains a non central abelian normal subgroup. As a consequence, we provide an extensive list of important families of groups which are generalized strongly monomial (Theorem \ref{t1}). In Theorem \ref{t2}, we provide an equivalent criterion for a generalized strong Shoda pair of $G$. The class $\mathcal{C}$ of groups as described above also appeared in our previous work \cite{BK}. We have noticed a misprint in the definition of $\mathcal{C}$ given in \cite{BK} where `quotient groups of $G$' should be read as `quotient groups of subgroups of $G$', i.e., the way class $\mathcal{C}$ is defined above. The theory developed in \cite{BK} is based on this property. We would like to point out that $\mathcal{C}$ is a large class of monomial groups as, in view of the work of Isaacs\cite{I} and Taketa\cite{T}, it seems difficult to construct example of monomial groups which are not in $\mathcal{C}$. \par In section 3, we generalize Proposition 3.4 of \cite{OdRS04} and provide an explicit description of the simple components of $ \mathbb{Q}G $ when $G$ is a generalized strongly monomial group (Theorem \ref{t3}). This also generalizes Theorem 5 of \cite{BK} on the simple components of $\mathbb{Q}G$ for $G \in \mathcal{C}.$ \par The knowledge of primitive central idempotents and simple components of $\mathbb{Q}G$ has a strong bearing upon understanding the unit group of the integral group ring $\mathbb{Z}G$. For details, we refer to \cite{JdR}, \cite{PS} and \cite{Ss}. In section 4, we prove that for a class of generalized strongly monomial groups, the group generated by generalized Bass units contain a subgroup of finite index in the group of central units of  $ \mathbb{Z}G$ (Theorem \ref{t4}). This generalizes the corresponding work on the central units of the integral group ring of a strongly monomial group by Jespers, Olteanu, del R{\'{\i}}o and Van Gelder (\cite{JOdRVG}, Theorem 5.1). Finally, we consider a class of groups in $\mathcal{C}$ containing abelian-by-supersolvable groups with the property that every cyclic subgroup of order not a divisor of 4 or 6 is subnormal in $G$. For a group $G$ in this class, we prove that the group generated by Bass units of $\mathbb{Z}G$ contains a subgroup of finite index in the group of central units of  $ \mathbb{Z}G$ (Theorem \ref{t5}), thus generalizing Theorem 3.2 of \cite{JOdRVG}. \section{Generalized strongly monomial groups} By $H\leq G$, $H \lneq G$, $H \unlhd G$, $H \unlhd$$\unlhd~G$, we mean, respectively, that $H$ is a subgroup, proper subgroup, normal subgroup, subnormal subgroup of $G$. Denote by \linebreak $[G:H]$, the index of $H$ in $G$. Also $N_{G}(H)$ denotes the normalizer of $H$ in $G$ and $\operatorname{core}_{G}(H)=\bigcap_{x \in G}x^{-1}Hx$ is the largest normal subgroup of $G$ contained in $H$. For $x,y\in G$, $[x,y]=x^{-1}y^{-1}xy$ is the commutator of $x$ and $y$, $x^{y}= y^{-1}xy$ and $\operatorname{Cen}_{G}(x)=\{g \in G~|~gx=xg\}$ is the centralizer of $x$ in $G$. Denote by $\operatorname{Irr}(G)$, the set of all complex irreducible characters of $G$. For $K\unlhd H\leq G$ with $H/K$ cyclic, $\operatorname{Lin}(H,K)$ denote the set of all linear characters on $H$ with kernel $K$. For a character $\chi$ of $G$, $\operatorname{ker}\chi =\{ g \in G~|~\chi(g)=\chi(1)\}$ and $\mathbb{Q}(\chi)$ denotes the field obtained by adjoining to $\mathbb{Q}$ the character values $\chi(g)$, $g\in G$. If $\psi$ is a character of a subgroup $H$ of $G$ and $x\in G$, then $\psi^{x}$ is the character of $H^{x}=x^{-1}Hx$ given by $\psi^{x}(g)=\psi(xgx^{-1})$, $g \in H^{x}$. Denote by $\psi^{G}$, the character $\psi$ induced to $G$. For a subgroup $A$ of $H$, $\psi_{A}$ denotes the restriction of $\psi$ to $A$. Let $I_{G}(\psi)=\{ g \in G~|~\psi^{g}=\psi\}$ be the inertia group of $\psi$ in $G$. For $\chi \in \operatorname{Irr}(G)$, $e_{\mathbb{Q}}(\chi)$ denotes the primitive central idempotent $\frac{\chi(1)}{|G|}\sum_{\sigma \in \operatorname{Gal}(\mathbb{Q}(\chi)/\mathbb{Q})} \sum_{g \in G}\sigma(\chi(g))g^{-1}$ of $\mathbb{Q}G$, where $\operatorname{Gal}(\mathbb{Q}(\chi)/\mathbb{Q})$ is the Galois group of $\mathbb{Q}(\chi)$ over $\mathbb{Q}$. For $\alpha \in \mathbb{Q}G$, $g \in G$ and an integer $m \geq 1$, $\alpha^{x}=x^{-1}\alpha x$ and $\alpha^{m}=\alpha\alpha \cdots \alpha$ ($m$ times).\par For $K\unlhd H\leq G$, define:$$\widehat{H}:=\frac{1}{|H|}\displaystyle\sum_{h \in H}h,$$ $$\varepsilon(H,K):=\left\{\begin{array}{ll}\widehat{K}, & \hbox{$H=K$;} \\\prod(\widehat{K}-\widehat{L}), & \hbox{otherwise,}\end{array}\right.$$ where $L$ runs over all the minimal normal subgroups of $H$ containing $K$ properly, and $$e(G,H,K):= {\rm~the~sum~of~all~the~distinct~}G{\rm {\tiny{\operatorname{-}}} conjugates~of~}\varepsilon(H,K).$$ \par We recall the notion of a Shoda pair,  strong Shoda pair and  strongly monomial group introduced in \cite{OdRS04}. A pair $(H,K)$ of subgroups of $G$ is called a \textit{Shoda pair} of $G$ if the following conditions hold:\begin{description}\item[(i)] $K\unlhd H$, $H/K$ is cyclic;\item[(ii)] if $g \in G$ and $[H,g]\cap H\subseteq K$, then $g \in H$.\end{description}\indent From Shoda's theorem (\cite{CR}, Corollary 45.4), it follows that if $(H,K)$ is a Shoda pair of $G$ and $\lambda \in \operatorname{Lin}(H,K)$, then $\lambda^{G}$ is irreducible. In this case, by Theorem 2.1 of \cite{OdRS04}, $e_{\mathbb{Q}}(\lambda^{G})=\alpha e(G,H,K)$ and it is called \textit{the primitive central idempotent of }$\mathbb{Q}G$ \textit{realized by the Shoda pair }$(H,K)$ of $G$, where $\alpha=\frac{[\operatorname{Cen}_{G}(\varepsilon(H,K)):H]}{[\mathbb{Q}(\lambda):\mathbb{Q}(\lambda^{G})]}$. Two Shoda pairs of $G$ are said to be \textit{equivalent} if they realize the same primitive central idempotent of $\mathbb{Q}G$. A set of representatives of distinct equivalence classes of Shoda pairs of $G$ is called a \textit{complete and irredundant set of Shoda pairs} of $G$. \par A Shoda pair $(H,K)$ of $G$ is called a \textit{strong Shoda pair} of $G$ if the following conditions hold:\begin{description}\item [(i)]$H \unlhd \operatorname{Cen}_{G}(\varepsilon(H,K))$;\item [(ii)] the distinct $G$-conjugates of $\varepsilon(H,K)$ are mutually orthogonal.\end{description}\par A group $G$ is called \textit{strongly monomial} if every irreducible character of $G$ is of the type $\lambda^{G}$, where $\lambda \in \operatorname{Lin}(H,K)$ and $(H,K)$ is a strong Shoda pair of $G$. For details on the theory of strong Shoda pairs and strongly monomial groups, see (\cite{JdR}, Chapter 3).\par An \textit{inductive source} of $G$, as defined by Dade \cite{D}, is an irreducible character $\psi$ of a subgroup of $G$ such that the induction to $G$ is an isometric bijection between the irreducible constituents of $\psi^{I_{G}(\psi)}$ and those of $\psi^{G}$. It is readily verified that if $\psi$ is an irreducible character of a subgroup of $G$ with $\psi^{G}$ irreducible then it is an inductive source of $G$. We call such an inductive source a \textit{direct inductive source} of $G$. By a \textit{strong inductive source} of $G$, we mean an irreducible character $\psi$ of a subgroup $S$ of $G$ which satisfies the following: \begin{description}\item [(i)] $\psi$ is a direct inductive source of $G$, i.e., $\psi^{G}$ is irreducible; \item [(ii)] $S \unlhd \operatorname{Cen}_{G}(e_{\mathbb{Q}}(\psi))$;\item [(iii)] the distinct $G$-conjugates of $e_{\mathbb{Q}}(\psi)$ are mutually orthogonal.\end{description}  It is not difficult to see that the  conditions (ii) and (iii) stated above are equivalent to saying that $e_{\mathbb{Q}}(\psi)$ and $e_{\mathbb{Q}}(\psi)^{g}$ are mutually orthogonal for every $g \in G\setminus N_{G}(S)$. \par  Let $(H,K)$ be a Shoda pair of $G$ and $\lambda \in \operatorname{Lin}(H,K)$. We say that $(H,K)$ is  a \textit{generalized strong Shoda pair} of $G$ if there is a chain $H=H_{0}\leq H_{1}\leq \cdots \leq H_{n}=G$ of subgroups of $G$ such that  $\lambda^{H_{i}}$ is a strong inductive source of $H_{i+1}$ for all $0\leq i \leq n-1$. We call such a chain of subgroups a \textit{strong inductive chain} from $H$ to $G$. Note that the definition of a generalized strong Shoda pair is independent of the choice of $\lambda$. Two generalized strong Shoda pairs of $G$ are called \textit{equivalent} if they are equivalent as Shoda pairs of $G$. A set of representatives of distinct equivalence classes of generalized strong Shoda pairs of $G$ is called a \textit{complete and irredundant set of generalized strong Shoda pairs} of $G$.\par A group $G$ is called \textit{generalized strongly monomial} if every irreducible character of $G$ is of the type $\lambda^{G}$, where $\lambda \in \operatorname{Lin}(H,K)$ and $(H,K)$ is a generalized strong Shoda pair of $G$. Note that if $G$ is a generalized strongly monomial group, then a complete and irredundant set of generalized strong Shoda pairs of $G$ is also a complete and irredundant set of Shoda pairs of $G$. \par Some observations are now in order. If $(H,K)$ is a strong Shoda pair of $G$ and $\lambda \in \operatorname{Lin}(H,K)$, then it follows immediately from the definition that $\lambda$ is a strong inductive source of $G$. Hence $(H,K)$ is a generalized strong Shoda pair of $G$ with $H\leq G$ a strong inductive chain. Consequently, the following statement holds: \begin{quote} All the strongly monomial groups are generalized strongly monomial.
\end{quote}
Next observe  that if $\psi \in \operatorname{Irr}(S)$ is a direct inductive source with $S \unlhd G$,  then it is a strong inductive source of $G$. This immediately yields the following: \begin{quote} A Shoda pair $(H,K)$ of $G$ with $H\unlhd$$\unlhd~G$ is a generalized strong Shoda pair of $G$. As a consequence, all the subnormally monomial groups are generalized strongly monomial.  \end{quote} It may be recalled that a group $G$ is called subnormally monomial if every complex irreducible character of $G$ is induced from a linear character of a subnormal subgroup of $G$. In (\cite{BK}, Proposition 1), we have shown that all subnormally monomial groups are contained in the class $\mathcal{C}$. We now show that this larger class $\mathcal{C}$ is contained in the class of generalized strongly monomial groups.\begin{theorem}\label{t1} The groups in class $\mathcal{C}$ are generalized strongly monomial. Consequently, the following groups are also generalized strongly monomial:\begin{description}\item [(i)] monomial Frobenius groups; \item [(ii)] Camina groups; \item [(iii)] groups such that all its non linear irreducible characters vanish only on the elements of order $p$, where $p$ is a fixed prime;\item [(iv)] solvable groups with the property that all its non linear irreducible characters of the same degree are Galois conjugate (this property holds in particular when all its non linear irreducible characters have distinct degree); \item [(v)] monomial groups whose non linear irreducible characters are induced from abelian subgroups; \item [(vi)] monomial groups such that all its non linear irreducible characters have same degree;\item [(vii)] monomial groups of odd order such that there are exactly two non linear irreducible characters of each degree; \item [(viii)] solvable groups with all its elements of prime power order; \item [(ix)] the class $\mathcal{A}$ of solvable groups with all its Sylow subgroups abelian; \item [(x)] $\mathcal{A}$-by-supersolvable groups, in particular, abelian-by-supersolvable groups; \item [(xi)] solvable groups $G$ satisfying the following condition: for all primes $p$ dividing the order of $G$ and for all subgroups $A$ of $G$, $O^{p}(A)$, the unique smallest normal subgroup of $A$ such that $A/O^{p}(A)$ is a $p$-group, has no central $p$-factor.\end{description}\end{theorem}\par We need some preparation to prove the theorem. A natural mechanism to produce a direct inductive source of $G$ is via Clifford theory. Let $\chi \in \operatorname{Irr}(G)$ and let $N \unlhd G$. Let $\varphi$ be an irreducible constituent of $\chi_{N}$. From Clifford's correspondence theorem (\cite{IM}, Theorem 6.11), there exist a unique $ \psi \in \operatorname{Irr}(I_{G}(\varphi))$ such that $\psi^{G} = \chi$. Such a direct inductive source $\psi$ is called the \textit{Clifford correspondent} of $\chi$ w.r.t. $\varphi$. We will show that if $\varphi$ is linear, then $\psi$ is a strong inductive source of $G$.\begin{prop}\label{p1} Let $\chi \in \operatorname{Irr}(G)$ and let $N \unlhd G$. Let $ \varphi \in \operatorname{Irr}(N)$ be a constituent of $\chi_{N}$ and let $\psi \in \operatorname{Irr}(I_{G}(\varphi))$ be the  Clifford correspondent of $\chi$ w.r.t. $\varphi$. If $\varphi$ is linear, then $\psi$ is a strong inductive source of $G$. \end{prop}\begin{lemma}\label{l1} Let $\chi \in \operatorname{Irr}(G)$ and let $N\unlhd G$. Suppose that $\chi_{N}=e\varphi$ for some \linebreak $\varphi \in \operatorname{Irr}(N)$ and a positive integer $e$. Then $$e_{\mathbb{Q}}(\chi)e_{\mathbb{Q}}(\varphi)=e_{\mathbb{Q}}(\chi)=e_{\mathbb{Q}}(\varphi)e_{\mathbb{Q}}(\chi).$$\end{lemma}{\bf Proof.} By Frobenius reciprocity, $\chi$ is an irreducible constituent of $\varphi^{G}$. Denote $\chi$ by $\chi_{1}$. Let $\chi_{1},~\chi_{2},\cdots,\chi_{m}$ be all the irreducible constituents of $\varphi^{G}$ so that  \begin{equation}\label{e1}\varphi^{G}=\sum_{1\leq i\leq m}\eta_{\chi_{i}}\chi_{i},~{\rm~where}~\eta_{\chi_{i}}\geq 1~{\rm are~integers}.\end{equation} Let $n$ be the exponent of $G$ and let $\xi$ be a complex primitive $n$th root of unity. Clearly $\mathbb{Q}(\chi_{i})\subseteq \mathbb{Q}(\xi)$ for all $i$. Define an equivalence relation on $\chi_{i},~1\leq i \leq m$, as follows: call $\chi_{i}$ to be equivalent to $\chi_{j}$, if there exist a $\sigma \in \operatorname{Gal}(\mathbb{Q}(\xi)/\mathbb{Q})$ such that $\sigma \circ \chi_{i}=\chi_{j}$. By reordering $\chi_{i}$'s if needed, we may assume that $\chi_{1},\cdots,~\chi_{k}$ form a set of representatives of distinct equivalence classes, where $k \leq m$. Denote by $\mathcal{G}_{i}$, the equivalence class of $\chi_{i}$, $1\leq i\leq k$. We can write eqn (\ref{e1}) as $$\varphi^{G}=\sum_{1\leq i \leq k}\sum_{\psi \in \mathcal{G}_{i}}\eta_{\psi}\psi.$$ Now,$$\begin{array}{lll}\displaystyle\sum_{\sigma \in \operatorname{Gal}(\mathbb{Q}(\xi)/\mathbb{Q})}\sum_{g \in G}\sigma \circ \varphi^{G}(g)g^{-1}&=&\displaystyle\sum_{g \in G}\sum_{1 \leq i\leq k}\sum_{\psi \in \mathcal{G}_{i}}\eta_{\psi}\sum_{\sigma \in \operatorname{Gal}(\mathbb{Q}(\xi)/\mathbb{Q})}\sigma\circ\psi(g)g^{-1}\vspace{.2cm}\\&=&\displaystyle\sum_{1 \leq i\leq k}\sum_{g \in G}t_{i}\sum_{\sigma \in \operatorname{Gal}(\mathbb{Q}(\xi)/\mathbb{Q})}\sigma\circ \chi_{i}(g)g^{-1}\vspace{.2cm}\\&=&\displaystyle\sum_{1\leq i\leq k}\frac{t_{i}[\mathbb{Q}(\xi):\mathbb{Q}(\chi_{i})]}{\chi_{i}(1)}|G|e_{\mathbb{Q}}(\chi_{i})\vspace{.2cm}\\
&=&\displaystyle\sum_{1\leq i\leq k}\alpha_{i}e_{\mathbb{Q}}(\chi_{i}),\end{array}$$ where $t_{i}=\sum_{\psi \in \mathcal{G}_{i}}\eta_{\psi}$ and $\alpha_{i} = \frac{t_{i}[\mathbb{Q}(\xi):\mathbb{Q}(\chi_{i})]|G|}{\chi_{i}(1)}.$\vspace{.2cm}\\ Also, $\varphi$ being $G$-invariant, we have
$$\begin{array}{lll}\displaystyle \sum_{\sigma \in \operatorname{Gal}(\mathbb{Q}(\xi)/\mathbb{Q})}\sum_{g \in G} \sigma\circ\varphi^{G}(g)g^{-1}&=&\displaystyle \sum_{\sigma \in \operatorname{Gal}(\mathbb{Q}(\xi)/\mathbb{Q})}\sum_{g \in N} [G:N]\sigma\circ\varphi(g)g^{-1}\vspace{.2cm}\\&=&\displaystyle\frac{|G|[\mathbb{Q}(\xi):\mathbb{Q}(\varphi)]}{\varphi(1)}e_{\mathbb{Q}}(\varphi)\vspace{.2cm}\\
&=&\displaystyle \beta e_{\mathbb{Q}}(\varphi),~{\rm where}~\beta=\frac{|G|[\mathbb{Q}(\xi):\mathbb{Q}(\varphi)]}{\varphi(1)}.\end{array}$$ Thus  \begin{equation}\label{e2} \beta e_{\mathbb{Q}}(\varphi)=\sum_{1\leq i \leq k}\alpha_{i}e_{\mathbb{Q}}(\chi_{i}).\end{equation} As $e_{\mathbb{Q}}(\chi_{i})$, $1\leq i \leq k$, are mutually orthogonal, multiplying eqn (\ref{e2}) by $e_{\mathbb{Q}}(\chi_{1})$, we have $$\beta e_{\mathbb{Q}}(\varphi)e_{\mathbb{Q}}(\chi_{1})=\alpha_{1}e_{\mathbb{Q}}(\chi_{1}).$$ Now both $e_{\mathbb{Q}}(\varphi)e_{\mathbb{Q}}(\chi_{1})$ and $e_{\mathbb{Q}}(\chi_{1})$ being idempotents, it follows immediately that $\alpha_{1}=\beta$ and $e_{\mathbb{Q}}(\varphi)e_{\mathbb{Q}}(\chi_{1})= e_{\mathbb{Q}}(\chi_{1})$. This finishes the proof.$~\Box$
 \vspace{.25cm}\\ {\bf Proof of Proposition \ref{p1}.} We will prove the result in steps:\vspace{.2cm}\\\underline{\textbf{Step I}} $I_{G}(\varphi)\unlhd N_{G}(\operatorname{ker}\varphi).$\\ Let $x \in N_{G}(\operatorname{ker}\varphi)$ and $t\in I_{G}(\varphi)$. Then $x^{-1}tx \in I_{G}(\varphi)$ if, and only if, $[x^{-1}tx,n]=x^{-1}[t,xnx^{-1}]x \in \operatorname{ker}\varphi$ for all $n \in N$. As $N\unlhd G$, the latter can be seen to hold.\vspace{.2cm}\\\underline{\textbf{Step II}} $e_{\mathbb{Q}}(\psi)e_{\mathbb{Q}}(\psi)^{g}=0$ if $g \in G\setminus N_{G}(\operatorname{ker}\varphi).$\\As $\psi \in \operatorname{Irr}(I_{G}(\varphi))$, it is easily seen using Clifford's theorem (\cite{IM}, Theorem 6.2) that $\psi_{N}=\psi(1)\varphi$. Hence by applying Lemma \ref{l1} to $\psi$, we obtain that $$e_{\mathbb{Q}}(\psi)e_{\mathbb{Q}}(\varphi)=e_{\mathbb{Q}}(\psi)=e_{\mathbb{Q}}(\varphi)e_{\mathbb{Q}}(\psi).$$ Consequently, for any $g \in G$, \begin{equation}\label{e3}e_{\mathbb{Q}}(\psi)e_{\mathbb{Q}}(\psi)^{g}=e_{\mathbb{Q}}(\psi)e_{\mathbb{Q}}(\varphi)e_{\mathbb{Q}}(\varphi)^{g}e_{\mathbb{Q}}(\psi)^{g}.
\end{equation} As $N\unlhd G$, $e_{\mathbb{Q}}(\varphi)$ and $e_{\mathbb{Q}}(\varphi)^{g}$ are both primitive central idempotents of $\mathbb{Q}N$. Since $g \not\in N_{G}(\operatorname{ker}\varphi)$, $e_{\mathbb{Q}}(\varphi)$ and $e_{\mathbb{Q}}(\varphi)^{g}$ can't be same and hence they are mutually orthogonal. Consequently, eqn (\ref{e3}) gives that $e_{\mathbb{Q}}(\psi)$ and $e_{\mathbb{Q}}(\psi)^{g}$ are mutually orthogonal.\vspace{.2cm}\\\underline{\textbf{Step III}} $\operatorname{Cen}_{G}(e_{\mathbb{Q}}(\psi))\leq N_{G}(\operatorname{ker}\varphi)$.\vspace{.15cm}\\ Let $x \in \operatorname{Cen}_{G}(e_{\mathbb{Q}}(\psi))$. Then $e_{\mathbb{Q}}(\psi)^{x}=e_{\mathbb{Q}}(\psi)$ and hence $e_{\mathbb{Q}}(\psi)^{x}$ and $e_{\mathbb{Q}}(\psi)$ can't be mutually orthogonal. Therefore, from the above step, it follows that $x\in N_{G}(\operatorname{ker}\varphi)$.\vspace{.2cm}\\\underline{\textbf{Step IV}} $e_{\mathbb{Q}}(\psi)e_{\mathbb{Q}}(\psi)^{g}=0$ if $g\in N_{G}(\operatorname{ker}\varphi)\setminus \operatorname{Cen}_{G}(e_{\mathbb{Q}}(\psi))$.\vspace{.15cm}\\ Let $g \in N_{G}(\operatorname{ker}\varphi)$. In view of step I, $e_{\mathbb{Q}}(\psi)$ and $e_{\mathbb{Q}}(\psi)^{g}$ are both primitive central idempotents of $\mathbb{Q}I_{G}(\varphi)$. Hence if $g \not\in \operatorname{Cen}_{G}(e_{\mathbb{Q}}(\psi))$, they must be mutually orthogonal.\par  From the above steps, it follows immediately that $I_{G}(\varphi)\unlhd \operatorname{Cen}_{G}(e_{\mathbb{Q}}(\psi))$ and the distinct $G$-conjugates of $e_{\mathbb{Q}}(\psi)$ are mutually orthogonal, i.e., $\psi$ is a strong inductive source of $G$. $~\Box$\begin{lemma}\label{l2} Let $\chi \in \operatorname{Irr}(G)$ be faithful and let $A$ be an abelian normal subgroup of $G$. Let $\varphi \in \operatorname{Irr}(A)$  be a constituent of $\chi_{A}$. Then $I_{G}(\varphi)$ is a proper subgroup of $G$ if, and only if, $A$ is non central in $G$.\end{lemma}{\bf Proof.} Let $\rho$ be the $\mathbb{C}$-representation of $G$ which affords the character $\chi$. As $A$ is an abelian normal subgroup of $G$, $I_{G}(\varphi)=G$ if, and only if, the restriction of $\rho$ to $A$ is isotypic, which is further equivalent to saying that $\rho(a)$ is homothety for each $a \in A$. Consequently, Lemma 2.27 of \cite{IM} gives that $I_{G}(\varphi)=G$ if, and only if, $A$ is central in $G$. This proves the desired result. $~\Box$\vspace{.25cm}\\ {\bf Proof of Theorem \ref{t1}.} Let $G \in \mathcal{C}$ and $\chi \in \operatorname{Irr}(G)$. We claim that $\chi = \lambda^{G},$ where $\lambda \in \operatorname{Lin}(H,K)$ and $(H,K)$ a generalized strong Shoda pair of $G$. If $\chi$ is linear, this is trivially true. Suppose $\chi$ is non linear. Let $A_{1}/\operatorname{ker}\chi$ be a non central abelian normal subgroup of $G/\operatorname{ker}\chi$. As $G \in \mathcal{C}$, such a subgroup exists. Let $\varphi_{1}$ be an irreducible constituent of $\chi_{A_{1}}$. Observe that $\varphi_{1}$ is linear as $A_{1}/\operatorname{ker}\chi$ is abelian and $\operatorname{ker}\chi\leq \operatorname{ker}\varphi_{1}$. Let $\chi_{0}=\chi$, $L_{0}=G$ and $L_{1}=I_{G}(\varphi_{1})$. By Lemma \ref{l2}, $L_{1}\lneq G$. Let $\chi_{1} \in \operatorname{Irr}(L_{1})$ be the Clifford correspondent of $\chi_{0}$ w.r.t. $\varphi_{1}$. Then $\chi_{0}=\chi_{1}^{L_{0}}$. Proposition \ref{p1} yields that $\chi_{1}$ is a strong inductive source of $G$. If $\chi_{1}$ is linear, stop. Otherwise, we repeat the process for $\chi_{1}$, which is possible because $L_{1} \in \mathcal{C}$. So we have $L_{2}\lneq L_{1}$, $\chi_{2} \in \operatorname{Irr}(L_{2})$ such that $\chi_{2}$ is a strong inductive source of $L_{1}$ and moreover $\chi_{2}^{L_{1}}=\chi_{1}$. If $\chi_{2}$ is non linear, continue. This process has to stop because we are getting a strictly descending chain of subgroups. Suppose the process stops at $n$th step. Then there are subgroups $G=L_{0}\gneq L_{1}\gneq \cdots \gneq L_{n}$ with $\chi_{i} \in \operatorname{Irr}(L_{i})$, $\chi_{n}$ linear, such that $\chi_{i}^{L_{i-1}}=\chi_{i-1}$, $\chi_{i}$ is a strong inductive source of $L_{i-1}$ for all $1\leq i\leq n$. These properties imply that  $(L_{n},\operatorname{ker}\chi_{n})$ is a generalized strong Shoda pair of $G$ and  $\chi = \chi_{n}^{G}$. Consequently, the claim  follows by taking $(H,K)=(L_{n}, \operatorname{ker}\chi_{n})$ and $\lambda = \chi_{n}$. We now proceed to prove the statements (i)-(xi). We will show that the groups of type (i)-(vii) are subnormally monomial and those of type (viii)-(xi) belong to $\mathcal{C}$.\\(i) This follows from (\cite{H}, Proposition 2.8).\vspace{.2cm}\\ (ii) If $G$ is a Camina group, then either of the following holds (see \cite{L}):\\(a) $G$ is a $p$-group;\\(b) $G$ is a Frobenius group whose complement is either cyclic or isomorphic to quaternion.\\If (a) holds, i.e., $G$ is a p-group, then every subgroup of $G$ is subnormal in $G$ and hence $G$ is subnormally monomial. If (b) holds, then again $G$ is subnormally monomial using (\cite{H}, Proposition 2.8).\vspace{.2cm}\\(iii) Suppose that every non linear irreducible character of $G$ vanishes only on the elements of order $p$. By (\cite{BDS}, Theorem B), either of the following holds:\\(a) $G$ is a $p$-group of exponent $p$;\\(b) $G=E\times F$, where $E$ is an elementary abelian $p$-group and $F$ is a Frobenius group with a Frobenius complement of order $p$.\\One can see that in both the cases, $G$ is subnormally monomial.\vspace{.2cm}\\(iv) If $G$ is a solvable group such that all its non linear irreducible characters of the same degree are Galois conjugate, then by (\cite{DY}, Theorem A), either $G$ is a $p$-group or $G$ is a Frobenius monomial group and hence the result follows.\vspace{.2cm}\\(v) Suppose $\chi$ is an irreducible character of $G$ induced from a linear character of an abelian subgroup, say $A$. Then Theorem A of \cite{R} yields that $A$ must be subnormal in $G$. Consequently, $G$ is subnormally monomial.\vspace{.2cm}\\(vi) From (\cite{gR}, Lemma 2.2), it follows that such a group is either a $p$-group or its every non linear irreducible character is induced from a linear character of an elementary abelian $p$-group. In both the cases, we are done using (v). \vspace{.2cm}\\(vii) This follows immediately from the corollary to Theorem 1 of \cite{CH}.\vspace{.2cm}\\(viii) Let $G$ be a solvable group such that all its elements are of prime power order. We first assert that either $G$ is abelian or there exists a non central abelian normal subgroup of $G$. If $G$ is a $p$-group, then this can be clearly seen to hold. We may assume that $G$ is not a $p$-group. Let $n$ be the derived length of $G$ and let $G^{(n-1)}$ be the $(n-1)$th derived subgroup of $G$. Clearly, $G^{(n-1)}$ is an abelian normal subgroup of $G$. We'll show that it is non central. This follows if we show that $\mathcal{Z}(G)$, the center of $G$, is trivial. Suppose $\mathcal{Z}(G)$ is non trivial.  Consider $g \in \mathcal{Z}(G)$ of prime order, say $p$. As $G$ is not a $p$-group, it has an element $h$ of prime order $q$, where $q\neq p$. Consequently, $gh$ has order $pq$, which is a contradiction to the hypothesis on $G$.  This proves the assertion.  Next, we see that  the hypothesis on $G$ is clearly satisfied by all the subquotients of $G$.  Hence the assertion implies that every subquotient of $G$ is either abelian or it contains a non central abelian normal subgroup of $G$. This proves that $G \in \mathcal{C}$.\vspace{.2cm}\\(ix), (x), (xi) follow from (\cite{BK}, Proposition 1).$~\Box$
\vspace{.25cm} \par The following result provides an equivalent criterion for a generalized strong Shoda pair. \begin{theorem}\label{t2} Let $K\leq H$ be subgroups of $G$. Then $(H,K)$ is a generalized strong Shoda pair of $G$ if, and only if, the following conditions hold: \begin{description}\item [(i)] $K\unlhd H$ with $H/K$ cyclic;\item [(ii)] there is a chain $H=H_{0}\leq H_{1}\leq \cdots \leq H_{n}=G$ of subgroups of $G$ such that for all $0 \leq i \leq n-1$:\begin{description} \item [(a)] $H_{i}\unlhd \operatorname{Cen}_{H_{i+1}}(e(H_{i},H,K))$; \item [(b)] the distinct $H_{i+1}$-conjugates of $e(H_{i},H,K)$ are mutually orthogonal;\item [(c)] $(\lambda^{H_{i}})^{x}=\lambda^{H_{i}}$ for $x \in \operatorname{Cen}_{H_{i+1}}(e(H_{i},H,K))$ implies that $x \in H_{i}$, where $\lambda$ is a linear character on $H$ with kernel $K$.\end{description}\end{description}\end{theorem} To prove the theorem, we need the following natural generalization of (\cite{OdRS04}, Theorem 2.1). \begin{lemma}\label{l3} If $\psi \in \operatorname{Irr}(S)$ is a direct inductive source of $G$, then
	$$e_{\mathbb{Q}}(\psi^{G})=\frac{[\operatorname{Cen}_{G}(e_{\mathbb{Q}}(\psi)):S]}{[\mathbb{Q}(\psi):\mathbb{Q}(\psi^{G})]}\sum_{g \in T}e_{\mathbb{Q}}(\psi)^{g},$$ where $T$ is a right transversal of $\operatorname{Cen}_{G}(e_{\mathbb{Q}}(\psi))$ in $G$. Moreover, if $\psi$ is a strong inductive source of $G$, then $e_{\mathbb{Q}}(\psi^{G})$ is the sum of all the distinct $G$-conjugates of $e_{\mathbb{Q}}(\psi)$.\end{lemma} {\bf Proof.} Denote by $\operatorname{Aut}(\mathbb{C})$, the group of automorphisms of the field $\mathbb{C}$ of complex numbers. Consider the action of $\operatorname{Aut}(\mathbb{C})$ and $G$ on $\mathbb{C}G$ respectively as follows:$$ \sigma.(\sum_{g \in G}a_{g}g)=\sum_{g \in G}\sigma(a_{g})g,~~\sigma \in \operatorname{Aut}(\mathbb{C}),~\sum_{g \in G}a_{g}g \in \mathbb{C}G;$$ $$x.g=g^{-1}xg,~~x \in \mathbb{C}G,~g \in G.$$ These actions satisfy $\sigma.(x.g)= (\sigma.x).g$ for all $\sigma \in \operatorname{Aut}(\mathbb{C}), x \in \mathbb{C}G, g \in G.$ Hence $\operatorname{Aut}(\mathbb{C})\times G$ acts on the left of $\mathbb{C}G$ by $(\sigma,g).x=\sigma.x.g^{-1}$, where $(\sigma,g) \in \operatorname{Aut}(\mathbb{C})\times G$ and $x \in \mathbb{C}G$. Also $G$ acts on the set of all irreducible characters of subgroups of $G$ by conjugation, i.e., $\chi.g = \chi^{g}$ for $\chi$ an irreducible character of a subgroup of $G$ and $g \in G$. Denote $e(\psi)$ by $e$. It can be easily seen that the $G$-stablizer of $e$ is same as the $G$-stablizer of $\psi$. However, as $\psi^{G}$ is irreducible, Mackey's irreducibility criterion yields that the $G$-stabilizer of $\psi$ is $S$. Therefore, the $G$-orbit of $e$ is $\{e.g~|~g \in T'\}$, where $T'$ is a right transversal of $S$ in $G$. In view of (\cite{Y}, Proposition 1.1), the $\operatorname{Aut}(\mathbb{C})$-orbit of $e$ is $\{\sigma.e~|~\sigma \in \operatorname{Gal}(\mathbb{Q}(\psi)/\mathbb{Q})\}$. Now consider the elements $\sigma.e.g$, where $\sigma \in \operatorname{Gal}(\mathbb{Q}(\psi)/\mathbb{Q})$ and $g \in T'$. Summing these elements in two different ways, as done in the proof of Theorem 2.1 of \cite{OdRS04}, it turns out that $$e_{\mathbb{Q}}(\psi^{G})=\frac{[\operatorname{Cen}_{G}(e_{\mathbb{Q}}(\psi)):S]}{[\mathbb{Q}(\psi):\mathbb{Q}(\psi^{G})]}\sum_{g \in T}e_{\mathbb{Q}}(\psi)^{g},$$ where $T$ is a right transversal of $\operatorname{Cen}_{G}(e_{\mathbb{Q}}(\psi))$ in $G$.  Finally, if $\psi$ is a strong inductive source of $G$, then the orthogonality of the distinct $G$-conjugates of $e_{\mathbb{Q}}(\psi)$ implies that $\sum_{g \in T}e_{\mathbb{Q}}(\psi)^{g}$ is an idempotent and consequently $\frac{[\operatorname{Cen}_{G}(e_{\mathbb{Q}}(\psi)):S]}{[\mathbb{Q}(\psi):\mathbb{Q}(\psi^{G})]}$ equals $1$, which completes the proof. $~\Box$\begin{lemma}\label{l4} If $(H,K)$ is a generalized strong Shoda pair of $G$ and $\lambda \in \operatorname{Lin}(H,K)$, then $$e_{\mathbb{Q}}(\lambda^{G})\varepsilon(H,K)=\varepsilon(H,K)=\varepsilon(H,K)e_{\mathbb{Q}}(\lambda^{G}).$$\end{lemma}{\bf Proof.} Let $H$=$H_{0}\leq H_{1}\leq \cdots \leq H_{n}$=$G$ be a strong inductive chain from $H$ to $G$. By Lemma \ref{l3}, we have $$e_{\mathbb{Q}}(\lambda^{H_{i+1}})e_{\mathbb{Q}}(\lambda^{H_{i}})=e_{\mathbb{Q}}(\lambda^{H_{i}})=e_{\mathbb{Q}}(\lambda^{H_{i}})
e_{\mathbb{Q}}(\lambda^{H_{i+1}}) {\rm~for~all~} 0\leq i \leq n-1.$$ Hence it follows that  $e_{\mathbb{Q}}(\lambda^{H_{i+1}})e_{\mathbb{Q}}(\lambda)=e_{\mathbb{Q}}(\lambda)=e_{\mathbb{Q}}(\lambda)e_{\mathbb{Q}}(\lambda^{H_{i+1}})$ for all $0\leq i \leq n-1.$ In particular, $e_{\mathbb{Q}}(\lambda^{G})e_{\mathbb{Q}}(\lambda)=e_{\mathbb{Q}}(\lambda)=e_{\mathbb{Q}}(\lambda)e_{\mathbb{Q}}(\lambda^{G})$. This yields the desired result as $e_{\mathbb{Q}}(\lambda)=\varepsilon(H,K)$. $~\Box$\vspace{.25cm}\\ {\bf Proof of Theorem \ref{t2}.} Suppose $(H,K)$ is a generalized strong Shoda pair of $G$. Using Mackey's irreducibility criterion, it can be readily verified that the `only if' conditions hold.\par Conversely, assume that the conditions (i) and (ii) hold. We'll  show by induction that $(H,K)$ is a generalized strong Shoda pair of $H_{i}$ for all $0\leq i\leq n$. For $i=0$, the result is trivially true. Assume that $(H,K)$ is a generalized strong Shoda pair of $H_{i}$ for some $0\leq i\leq n-1$. To show that $(H,K)$ is a generalized strong Shoda pair of $H_{i+1}$, in view of conditions (ii)(a) and (ii)(b), it is enough to see that $\lambda^{H_{i+1}}$ is irreducible, equivalently, the following hold: \begin{equation}\label{e4}\lambda_{H\cap H^{x}}\neq \lambda^{x}_{H\cap H^{x}}~~{\rm for~all}~x \in H_{i+1}\setminus H.\end{equation} Suppose $x \in H_{i+1}\setminus H.$ The following two cases arise:\vspace{.2cm}\\\underline{\textbf{Case I}}~~$x \in \operatorname{Cen}_{H_{i+1}}(e_{\mathbb{Q}}(\lambda^{H_{i}}))\setminus H$.\\ By induction, $\lambda^{H_{i}}$ is irreducible. In view of (ii)(a) and (ii)(c), it follows using Mackey's irreducibility criterion that $\psi^{\operatorname{Cen}_{H_{i+1}}(e_{\mathbb{Q}}(\lambda^{H_{i}}))}$ is irreducible. Consequently, eqn (\ref{e4}) holds for all $x \in \operatorname{Cen}_{H_{i+1}}(e_{\mathbb{Q}}(\lambda^{H_{i}}))\setminus H$.\vspace{.2cm}\\\underline{\textbf{Case II}} $x \in H_{i+1}\setminus \operatorname{Cen}_{H_{i+1}}(e_{\mathbb{Q}}(\lambda^{H_{i}}))$.\\ Suppose $\lambda_{H\cap H^{x}}=\lambda^{x}_{H\cap H^{x}}$. We have $$e(\lambda)e(\lambda)^{x}=\frac{1}{|H|^{2}}\sum_{h_{1},h_{2}\in H}\lambda(h_{1}^{-1})\lambda(h_{2}^{-1})h_{1}x^{-1}h_{2}x.$$ Thus the coefficient of $1$ in $e(\lambda)e(\lambda)^{x}$ equals $$\begin{array}{lll} \displaystyle \frac{1}{|H|^{2}}\sum_{h\in H\cap H^{x}}\lambda(h^{-1})\lambda(xhx^{-1})&=&\displaystyle \frac{1}{|H|^{2}}\sum_{h\in H\cap H^{x}}\lambda(h^{-1})\lambda(h)~({\rm as}~ \lambda_{H\cap H^{x}} = \lambda^{x}_{H \cap H^{x}})\\&=&\displaystyle\frac{1}{|H|^{2}}\sum_{h \in H\cap H^{x}}\lambda(1)\\&=&\displaystyle\frac{1}{|H|^{2}}|H\cap H^{x}|\neq 0.\end{array}$$ Therefore, $e(\lambda)e(\lambda)^{x}\neq 0.$ As $e(\lambda)\varepsilon(H,K)=e(\lambda)=\varepsilon(H,K)e(\lambda)$, it immediately follows that $$\varepsilon(H,K)\varepsilon(H,K)^{x}\neq 0.$$ Since $e_{\mathbb{Q}}(\lambda^{H_{i}})\varepsilon(H,K)=\varepsilon(H,K)=\varepsilon(H,K)e_{\mathbb{Q}}(\lambda^{H_{i}})$ (by Lemma \ref{l4}), we obtain that $e_{\mathbb{Q}}(\lambda^{H_{i}})e_{\mathbb{Q}}(\lambda^{H_{i}})^{x}\neq 0.$ This is a contradiction to (ii)(b) and completes the proof of the theorem. $~\Box$ \section{Simple components} For a ring $R$, let $\mathcal{U}(R)$ be the unit group of $R$ and $M_{n}(R)$ the ring of $n\times n$ matrices over $R$. Denote by $R*_{\tau}^{\sigma}G$, the crossed product of the group $G$ over the ring $R$ with action $\sigma$ and twisting $\tau$ (see \cite{JdR}, Chapter 2).\par Let $S$ be a subgroup of $G$ and let $\psi \in \operatorname{Irr}(S)$ be a strong inductive source of $G$. For notational convenience, denote $\operatorname{Cen}_{G}(e_{\mathbb{Q}}(\psi))$ by $C$. For each $x \in C/S$, let $\overline{x}\in C$ be a fixed inverse image of $x$ under the natural map $C\rightarrow C/S.$ Since $\overline{x}$ centralizes $e_{\mathbb{Q}}(\psi)$, we have $\overline{x}\mathbb{Q}Se_{\mathbb{Q}}(\psi)\overline{x}^{-1}=\mathbb{Q}Se_{\mathbb{Q}}(\psi)$. Consequently, there is a map \linebreak $\sigma_{S}:C/S\rightarrow \operatorname{Aut}(\mathbb{Q}Se_{\mathbb{Q}}(\psi))$, which maps $x$ to the conjugation automorphism $(\sigma_{S})_{x}$ on $\mathbb{Q}Se_{\mathbb{Q}}(\psi)$ induced by $\overline{x}$. Also $\overline{x}r=(\overline{x}r\overline{x}^{-1})\overline{x}=(\sigma_{S})_{x}(r)\overline{x}$ for $x \in C/S$ and $r \in \mathbb{Q}Se_{\mathbb{Q}}(\psi)$, and we see that an action occurs. Furthermore for $x,y \in C/S$, we have $S\overline{x}S\overline{y}=S\overline{xy}$ and so $\overline{x}.\overline{y}=s\overline{xy},$ where $s \in S$. Thus we have a map $\tau_{S}:
 C/S\times C/S\rightarrow \mathcal{U}(\mathbb{Q}Se_{\mathbb{Q}}(\psi))$ such that $\tau_{S}(x,y) = se_{\mathbb{Q}}(\psi) \in \mathcal{U}(\mathbb{Q}Se_{\mathbb{Q}}(\psi)).$ Finally we can check that $\mathbb{Q}Ce_{\mathbb{Q}}(\psi)= \mathbb{Q}Se_{\mathbb{Q}}(\psi)\ast^{\sigma_{S}}_{\tau_{S}}C/S.$ \begin{prop}\label{p2} Let $\psi \in \operatorname{Irr}(S)$ be a strong inductive source of $G$ and let \linebreak $T=\{t_{i}~|~1 \leq i\leq m \}$ be a right transversal of $C$ in $G$, where $C=\operatorname{Cen}_{G}(e_{\mathbb{Q}}(\psi))$. Then the following hold:\begin{description}\item [(i)] $\mathbb{Q}Ge_{\mathbb{Q}}(\psi^{G})$ is isomorphic to $M_{m}(\mathbb{Q}Ce_{\mathbb{Q}}(\psi))$ and the map which defines this isomorphism is given by $\alpha \mapsto (\alpha_{ij})_{m\times m},$ where $\alpha_{ij}=e_{\mathbb{Q}}(\psi)t_{j}\alpha t_{i}^{-1}e_{\mathbb{Q}}(\psi)$; \item [(ii)] the map $\alpha \mapsto \sum_{t \in T}\alpha^{t}$ defines an isomorphism from $\mathcal{Z}(\mathbb{Q}Ce_{\mathbb{Q}}(\psi))$ to \linebreak $\mathcal{Z}(\mathbb{Q}Ge_{\mathbb{Q}}(\psi^{G}))$; \item [(iii)] $C/S$ acts on $\mathcal{Z}(\mathbb{Q}Se_{\mathbb{Q}}(\psi))$ by conjugation and $\mathcal{Z}(\mathbb{Q}Se_{\mathbb{Q}}(\psi))^{C/S} \linebreak = \{ r \in \mathcal{Z}(\mathbb{Q}Se_{\mathbb{Q}}(\psi))~|~(\sigma_{S})_{x}(r)=r~\forall x \in C/S\}$is a subring of $\mathcal{Z}(\mathbb{Q}Ce_{\mathbb{Q}}(\psi))$.\end{description}\end{prop}{\bf Proof.} (i) For $\alpha \in \mathbb{Q}Ge_{\mathbb{Q}}(\psi^{G})$, let $T_{\alpha}$ be the $\mathbb{Q}G$ endomorphism of $\mathbb{Q}Ge_{\mathbb{Q}}(\psi^{G})$ given by $\beta \mapsto \beta\alpha$. One can easily check that $\vartheta_{1}:~\mathbb{Q}Ge_{\mathbb{Q}}(\psi^{G})\rightarrow \operatorname{End}_{\mathbb{Q}G}(\mathbb{Q}Ge_{\mathbb{Q}}(\psi^{G}))$ given by $\alpha \mapsto T_{\alpha}$ is a ring isomorphism. We have $\displaystyle\mathbb{Q}Ge_{\mathbb{Q}}(\psi^{G})=\oplus_{k=1}^{m}\mathbb{Q}Ge_{\mathbb{Q}}(\psi)^{t_{k}}$. Moreover, for all $1\leq k\leq m$, $\mathbb{Q}Ge_{\mathbb{Q}}(\psi)^{t_{k}}\cong \mathbb{Q}Ge_{\mathbb{Q}}(\psi)$ as $\mathbb{Q}G$ modules, and the map which defines this isomorphism is given by $\gamma \mapsto \gamma t_{k}^{-1}$. As a consequence, it follows that  $\mathbb{Q}Ge_{\mathbb{Q}}(\psi^{G})$ is isomorphic to direct sum of $m$ copies of $\mathbb{Q}Ge_{\mathbb{Q}}(\psi)$ as $\mathbb{Q}G$ modules. Let $\pi_{i}:~\mathbb{Q}Ge_{\mathbb{Q}}(\psi^{G})\rightarrow \mathbb{Q}Ge_{\mathbb{Q}}(\psi)$ be the projection of $\mathbb{Q}Ge_{\mathbb{Q}}(\psi^{G})$ onto its $i$th component and let $\varepsilon_{j}:~\mathbb{Q}Ge_{\mathbb{Q}}(\psi)\rightarrow \mathbb{Q}Ge_{\mathbb{Q}}(\psi^{G})$ be the inclusion of $\mathbb{Q}Ge_{\mathbb{Q}}(\psi)$ into $j$th component, where $1\leq i,j\leq m$. By Lemma 2.6.11 of \cite{Ss}, it follows that the map $\vartheta_{2}:~f\mapsto (\pi_{i}\circ f\circ \varepsilon_{j})_{m \times m}$ is a ring isomorphism from $\operatorname{End}_{\mathbb{Q}G}(\mathbb{Q}Ge_{\mathbb{Q}}(\psi^{G}))$ to $M_{m}(\operatorname{End}_{\mathbb{Q}G}(\mathbb{Q}Ge_{\mathbb{Q}}(\psi)))$. One further checks that $\vartheta_{3}:~M_{m}(\operatorname{End}_{\mathbb{Q}G}(\mathbb{Q}Ge_{\mathbb{Q}}(\psi)))\rightarrow M_{m}(e_{\mathbb{Q}}(\psi)\mathbb{Q}Ge_{\mathbb{Q}}(\psi))$ given by $M_{m}(f_{ij})\mapsto M_{m}(f_{ij}(e_{\mathbb{Q}}(\psi)))$ is a ring isomorphism. As the distinct $G$-conjugates of $e_{\mathbb{Q}}(\psi)$ are mutually orthogonal, we have $e_{\mathbb{Q}}(\psi)\mathbb{Q}Ge_{\mathbb{Q}}(\psi)= e_{\mathbb{Q}}(\psi)\mathbb{Q}Ce_{\mathbb{Q}}(\psi)=\mathbb{Q}Ce_{\mathbb{Q}}(\psi)$. Consequently, $\vartheta=\vartheta_{3}\circ \vartheta_{2} \circ \vartheta_{1}$ is the required isomorphism from $\mathbb{Q}Ge_{\mathbb{Q}}(\psi^{G})$ to $M_{m}(\mathbb{Q}Ce_{\mathbb{Q}}(\psi))$. Little computation reveal that $\vartheta(\alpha) = (\alpha_{ij})_{m \times m}$, where $\alpha_{ij}=e_{\mathbb{Q}}(\psi)t_{j}\alpha t_{i}^{-1}e_{\mathbb{Q}}(\psi)$. This finishes the proof of (i).\vspace{.2cm}\\(ii) For $\alpha \in \mathcal{Z}(\mathbb{Q}Ce_{\mathbb{Q}}(\psi))$, one checks that the isomorphism $\vartheta$ given in (i) maps $\sum_{t \in T}\alpha^{t}$ to the scalar matrix with $\alpha$ on the diagonal. This proves (ii).\vspace{.2cm}\\(iii) This is easy to see. $~\Box$\vspace{.25cm}\para Let $(H,K)$ be a generalized strong Shoda pair of $G$ and $\lambda \in \operatorname{Lin}(H,K)$. Let $H=H_{0}\leq H_{1}\leq \cdots \leq H_{n}=G$ be a strong inductive chain from $H$ to $G$. We will recursively apply Proposition \ref{p2} to describe the simple component $\mathbb{Q}Ge_{\mathbb{Q}}(\lambda^{G})$. Set $C_{i}=\operatorname{Cen}_{H_{i+1}}(e_{\mathbb{Q}}(\lambda^{H_{i}})),~0\leq i\leq n-1$. As $\lambda$ is a strong inductive source of $H_{1}$, Proposition \ref{p2} gives the isomorphism $$\mathbb{Q}H_{1}e_{\mathbb{Q}}(\lambda^{H_{1}}) \cong M_{k_{1}}(\mathbb{Q}He_{\mathbb{Q}}(\lambda)\ast^{\sigma_{H_{0}}}_{\tau_{H_{0}}} C_{0}/H_{0}),$$  where $k_{1}=[H_{1}:C_{0}].$ The action and twisting $\sigma_{H_{1}}$ and $\tau_{H_{1}}$ on $\mathbb{Q}H_{1}e_{\mathbb{Q}}(\lambda^{H_{1}})$ induces the corresponding action and twisting on $M_{k_{1}}(\mathbb{Q}He_{\mathbb{Q}}(\lambda)\ast^{\sigma_{H_{0}}}_{\tau_{H_{0}}} C_{0}/H_{0})$, and again using Proposition \ref{p2}, we have the isomorphism $$\begin{array}{lll}\mathbb{Q}H_{2}e_{\mathbb{Q}}(\lambda^{H_{2}}) &\cong & M_{k_{2}}(\mathbb{Q}H_{1}e_{\mathbb{Q}}(\lambda^{H_{1}})\ast^{\sigma_{H_{1}}}_{\tau_{H_{1}}} C_{1}/H_{1})\\&\cong & M_{k_{2}}(M_{k_{1}}(\mathbb{Q}He_{\mathbb{Q}}(\lambda)\ast^{\sigma_{H_{0}}}_{\tau_{H_{0}}} C_{0}/H_{0})\ast^{\sigma_{H_{1}}}_{\tau_{H_{1}}} C_{1}/H_{1}),\end{array}$$ where $k_{2}=[H_{2}:C_{1}].$ Continue this process, the final step gives that {\footnotesize
 $$\mathbb{Q}H_{n}e_{\mathbb{Q}}(\lambda^{H_{n}}) \cong  M_{k_{n}}(M_{k_{n-1}}\cdots(M_{k_{1}}(\mathbb{Q}He_{\mathbb{Q}}(\lambda)\ast^{\sigma_{H_{0}}}_{\tau_{H_{0}}} C_{0}/H_{0})\ast^{\sigma_{H_{1}}}_{\tau_{H_{1}}}C_{1}/H_{1})\ast^{\sigma_{H_{2}}}_{\tau_{H_{2}}}
\cdots\ast^{\sigma_{H_{n-1}}}_{\tau_{H_{n-1}}}C_{n-1}/H_{n-1}).$$} Hence, with the forgoing notation, we have proved the following: \begin{theorem}\label{t3} Let $(H,K)$ be a generalized strong Shoda pair of $G$ and $\lambda \in \operatorname{Lin}(H,K)$.\vspace{.2cm}\\ Let $H=H_{0}\leq H_{1}\leq \cdots \leq H_{n}=G$ be a strong inductive chain from $H$ to $G$, then {\footnotesize $$\mathbb{Q}Ge_{\mathbb{Q}}(\lambda^{G}) \cong  M_{k_{n}}(M_{k_{n-1}}\cdots(M_{k_{1}}(\mathbb{Q}He_{\mathbb{Q}}(\lambda)\ast^{\sigma_{H_{0}}}_{\tau_{H_{0}}}C_{0}/H_{0})\ast^{\sigma_{H_{1}}}_{\tau_{H_{1}}}C_{1}
/H_{1})\ast^{\sigma_{H_{2}}}_{\tau_{H_{2}}}
\cdots\ast^{\sigma_{H_{n-1}}}_{\tau_{H_{n-1}}}C_{n-1}/H_{n-1}),$$} where $C_{i},~\sigma_{H_{i}},~\tau_{H_{i}},~k_{i}$ are as defined above. \end{theorem} \section{Central units of $\mathbb{Z}G$} We continue to use the notation developed in the previous section. The purpose of this section is to generalize the main work in \cite{JOdRVG}. We begin by generalizing an idea which originates from the proof of (\cite{JOdRVG}, Theorem 3.2). \par In order to proceed, we need to introduce a  terminology. Let  $ \psi \in \operatorname{Irr}(S)$ be a strong inductive source of $G$. Denote $\operatorname{Cen}_{G}(e_{\mathbb{Q}}(\psi))$ by $C$. In section 3, we have seen that $\mathbb{Q}Ce_{\mathbb{Q}}(\psi) = \mathbb{Q}Se_{\mathbb{Q}}(\psi)\ast^{\sigma_{S}}_{\tau_{S}}C/S.$ By Proposition \ref{p2}, $\mathcal{Z}(\mathbb{Q}Se_{\mathbb{Q}}(\psi))^{C/S}$ is a subring of $\mathcal{Z}(\mathbb{Q}Ce_{\mathbb{Q}}(\psi))$. Suppose $\psi$ is such that these are equal, i.e., $\mathcal{Z}(\mathbb{Q}Ce_{\mathbb{Q}}(\psi))=
\mathcal{Z}(\mathbb{Q}Se_{\mathbb{Q}}(\psi))^{C/S}$, then we say that $\psi$ \textit{has the least center property w.r.t.} $G$. A sufficient condition for $\psi$ to have the least center property w.r.t. $G$ is when $\sigma_{S}$ is such that $(\sigma_{S})_{x}$ is not inner in $\mathbb{Q}Se_{\mathbb{Q}}(\psi)$ for every non trivial $x \in C/S$. This is an immediate consequence of (\cite{JdR}, Lemma 2.6.1). Furthermore, the condition that $\sigma_{S}$ is such that $(\sigma_{S})_{x}$ is not inner in $\mathbb{Q}Se_{\mathbb{Q}}(\psi)$ for every non trivial $x \in C/S$  trivially holds when $\sigma_{S}$ is one-one and $\mathbb{Q}Se_{\mathbb{Q}}(\psi)$ is commutative.  As an example to illustrate this situation, consider a strong Shoda pair $(H, K)$ of $G$ and a linear character $\psi$ on $H$ with kernel $K$. Then $\psi$ is a strong inductive source of $G$ and it has the property that  $\sigma_{H}$ is one-one and $\mathbb{Q}He_{\mathbb{Q}}(\psi)= \mathbb{Q}H\varepsilon(H,K)$ is commutative, which yields that this $\psi$ has the least center property w.r.t. $G$.\par In the next lemma, we will show that if $\psi \in \operatorname{Irr}(S)$ is a strong inductive source of $G$ having the least center property w.r.t. $G$, then a subgroup of $\mathcal{Z}(\mathcal{U}(\mathbb{Z}S))\cap \mathcal{Z}(\mathcal{U}(\mathbb{Z}(1-e_{\mathbb{Q}}(\psi))+\mathbb{Z}Se_{\mathbb{Q}}(\psi)))$ which is of finite index in $\mathcal{Z}(\mathcal{U}(\mathbb{Z}(1-e_{\mathbb{Q}}(\psi))+\mathbb{Z}Se_{\mathbb{Q}}(\psi)))$ can be used to construct a subgroup of $\mathcal{Z}(\mathcal{U}(\mathbb{Z}G))\cap \mathcal{Z}(\mathcal{U}(\mathbb{Z}(1-e_{\mathbb{Q}}(\psi^{G}))+\mathbb{Z}Ge_{\mathbb{Q}}(\psi^{G})))$ which is of finite index in $\mathcal{Z}(\mathcal{U}(\mathbb{Z}(1-e_{\mathbb{Q}}(\psi^{G}))+\mathbb{Z}Ge_{\mathbb{Q}}(\psi^{G})))$.\begin{lemma}\label{l5} Let $\psi \in \operatorname{Irr}(S)$ be a strong inductive source of $G$ having the least center property w.r.t. $G$ and let $T$ be a right transversal of $\operatorname{Cen}_{G}(e_{\mathbb{Q}}(\psi))$ in $G$. If $A$ is a subgroup of $\mathcal{Z}(\mathcal{U}(\mathbb{Z}S))\cap \mathcal{Z}(\mathcal{U}(\mathbb{Z}(1-e_{\mathbb{Q}}(\psi))+\mathbb{Z}Se_{\mathbb{Q}}(\psi)))$ which is of finite index in $\mathcal{Z}(\mathcal{U}(\mathbb{Z}(1-e_{\mathbb{Q}}(\psi))+\mathbb{Z}Se_{\mathbb{Q}}(\psi)))$, then the elements $$\prod_{t \in T}\Big (\prod_{c \in \operatorname{Cen}_{G}(e_{\mathbb{Q}}(\psi))}u^{c}\Big)^{t},~~ u \in A,$$ form a subgroup of $\mathcal{Z}(\mathcal{U}(\mathbb{Z}G))\cap \mathcal{Z}(\mathcal{U}(\mathbb{Z}(1-e_{\mathbb{Q}}(\psi^{G}))+\mathbb{Z}Ge_{\mathbb{Q}}(\psi^{G})))$ and its index in $\mathcal{Z}(\mathcal{U}(\mathbb{Z}(1-e_{\mathbb{Q}}(\psi^{G}))+\mathbb{Z}Ge_{\mathbb{Q}}(\psi^{G})))$ is finite.
\end{lemma}{\bf Proof.} \underline{\textbf{Step I}} $\mathfrak{B}=\{ \prod_{c \in \operatorname{Cen}_{G}(e_{\mathbb{Q}}(\psi))}u^{c}~|~u \in A\}$ is a subgroup of finite index in  $\mathcal{U}(\mathbb{Z}(1-e_{\mathbb{Q}}(\psi))+\mathcal{Z}(\mathbb{Z}Se_{\mathbb{Q}}(\psi))^{\operatorname{Cen}_{G}(e_{\mathbb{Q}}(\psi))/S})$. \vspace{.2cm}\\ For an arbitrary $u \in \mathcal{Z}(\mathcal{U}(\mathbb{Z}(1-e_{\mathbb{Q}}(\psi))+\mathbb{Z}Se_{\mathbb{Q}}(\psi)))$, observe that $\prod_{c \in \operatorname{Cen}_{G}(e_{\mathbb{Q}}(\psi))}u^{c}$ is independent of the order of entries in the product. Moreover, $\prod_{c \in \operatorname{Cen}_{G}(e_{\mathbb{Q}}(\psi))}u^{c}$ belongs to $\mathcal{U}(\mathbb{Z}(1-e_{\mathbb{Q}}(\psi))+\mathcal{Z}(\mathbb{Z}Se_{\mathbb{Q}}(\psi))^{\operatorname{Cen}_{G}(e_{\mathbb{Q}}(\psi))/S}).$ With these observations, it is easy to see that $\mathfrak{B}$ is a subgroup of $\mathcal{U}(\mathbb{Z}(1-e_{\mathbb{Q}}(\psi))+\mathcal{Z}(\mathbb{Z}Se_{\mathbb{Q}}(\psi))^{\operatorname{Cen}_{G}(e_{\mathbb{Q}}(\psi))/S})$. To see that the index is finite, it is enough to show that a power of any element of  $\mathcal{U}(\mathbb{Z}(1-e_{\mathbb{Q}}(\psi))+\mathcal{Z}(\mathbb{Z}Se_{\mathbb{Q}}(\psi))^{\operatorname{Cen}_{G}(e_{\mathbb{Q}}(\psi))/S})$ belongs to $\mathfrak{B}$, as $\mathcal{U}(\mathbb{Z}(1-e_{\mathbb{Q}}(\psi))+\mathcal{Z}(\mathbb{Z}Se_{\mathbb{Q}}(\psi))^{\operatorname{Cen}_{G}(e_{\mathbb{Q}}(\psi))/S})$ is  a subgroup of the finitely generated abelian group $\mathcal{Z}(\mathcal{U}(\mathbb{Z}(1-e_{\mathbb{Q}}(\psi))+\mathbb{Z}Se_{\mathbb{Q}}(\psi)))$. Consider an arbitrary $u \in \mathcal{U}(\mathbb{Z}(1-e_{\mathbb{Q}}(\psi))+\mathcal{Z}(\mathbb{Z}Se_{\mathbb{Q}}(\psi))^{\operatorname{Cen}_{G}(e_{\mathbb{Q}}(\psi))/S})$. We can write  $u= \pm (1-e_{\mathbb{Q}}(\psi))+u'$, where $u' \in \mathcal{U}(\mathcal{Z}(\mathbb{Z}Se_{\mathbb{Q}}(\psi))^{\operatorname{Cen}_{G}(e_{\mathbb{Q}}(\psi))/S})$. Suppose $u = 1-e_{\mathbb{Q}}(\psi)+u'$. Clearly, $$\prod_{c \in \operatorname{Cen}_{G}(e_{\mathbb{Q}}(\psi))}u^{c}=1-e_{\mathbb{Q}}(\psi)+\prod_{c \in \operatorname{Cen}_{G}(e_{\mathbb{Q}}(\psi))}u'^{c}=1-e_{\mathbb{Q}}(\psi)+\prod_{c \in \operatorname{Cen}_{G}(e_{\mathbb{Q}}(\psi))}u'=u^{t},$$ where $t=|\operatorname{Cen}_{G}(e_{\mathbb{Q}}(\psi))|$. Also, if $k$ is the index of $A$ in $\mathcal{Z}(\mathcal{U}(\mathbb{Z}(1-e_{\mathbb{Q}}(\psi))+\mathbb{Z}Se_{\mathbb{Q}}(\psi)))$, then $(\prod_{c \in \operatorname{Cen}_{G}(e_{\mathbb{Q}}(\psi))}u^{c})^{k}= \prod_{c \in \operatorname{Cen}_{G}(e_{\mathbb{Q}}(\psi))}(u^{k})^{c}$ belongs to $\mathfrak{B}$. Hence $u^{tk} \in \mathfrak{B}$, as desired. If $u = -(1-e_{\mathbb{Q}}(\psi))+u'$, then $u^{2}=1-e_{\mathbb{Q}}(\psi)+u'^{2}$ and hence $u^{2tk} \in \mathfrak{B}$. This proves step I. \par Let $T$ be a right transversal of $\operatorname{Cen}_{G}(e_{\mathbb{Q}}(\psi))$ in $G$. Our next goal is the following: \vspace{.2cm} \\ \underline{\textbf{Step II}} $\mathfrak{C}=\{ \prod_{t \in T}v^{t}~|~v \in \mathfrak{B}\}$ is a subgroup of $\mathcal{Z}(\mathcal{U}(\mathbb{Z}G))\cap \mathcal{Z}(\mathcal{U}(\mathbb{Z}(1-e_{\mathbb{Q}}(\psi^{G}))+\mathbb{Z}Ge_{\mathbb{Q}}(\psi^{G})))$ and its index in $\mathcal{Z}(\mathcal{U}(\mathbb{Z}(1-e_{\mathbb{Q}}(\psi^{G}))+\mathbb{Z}Ge_{\mathbb{Q}}(\psi^{G})))$ is finite. \vspace{.2cm} \\ Consider an arbitrary $v \in \mathfrak{B}$. We can write $v =\pm (1-e_{\mathbb{Q}}(\psi))+v'$, where $v' \in \mathcal{U}(\mathcal{Z}(\mathbb{Z}Se_{\mathbb{Q}}(\psi))^{\operatorname{Cen}_{G}(e_{\mathbb{Q}}(\psi))/S})$. Suppose $v = 1-e_{\mathbb{Q}}(\psi)+v'$. As $\psi$ is a strong inductive source of $G$, the distinct $G$-conjugates of $e_{\mathbb{Q}}(\psi)$ are mutually orthogonal and moreover, by Lemma \ref{l3}, $\sum_{t \in T}e_{\mathbb{Q}}(\psi)^{t}=e_{\mathbb{Q}}(\psi^{G})$. Hence $\prod_{t \in T}v^{t}=1-\sum_{t \in T}e_{\mathbb{Q}}(\psi)^{t}+\sum_{t \in T}v'^{t}= 1-e_{\mathbb{Q}}(\psi^{G})+\sum_{t \in T}v'^{t}$. This shows that $\prod_{t \in T}v^{t} $  is independent of the order of entries in the product. Moreover, this product is also independent of the transversal $T$. Further, by Proposition \ref{p2}, it follows that $\sum_{t \in T}v'^{t} \in \mathcal{Z}(\mathcal{U}(\mathbb{Z}Ge_{\mathbb{Q}}(\psi^{G})))$. Consequently, $\prod_{t \in T}v^{t} \in \mathcal{Z}(\mathcal{U}(\mathbb{Z}(1-e_{\mathbb{Q}}(\psi^{G}))+\mathbb{Z}Ge_{\mathbb{Q}}(\psi^{G})))$. Similarly, if $v = -(1-e_{\mathbb{Q}}(\psi))+v'$, then it can be verified that $\prod_{t \in T}v^{t} \in \mathcal{Z}(\mathcal{U}(\mathbb{Z}(1-e_{\mathbb{Q}}(\psi^{G}))+\mathbb{Z}Ge_{\mathbb{Q}}(\psi^{G})))$. Hence $\mathfrak{C}$ is contained in $\mathcal{Z}(\mathcal{U}(\mathbb{Z}(1-e_{\mathbb{Q}}(\psi^{G}))+\mathbb{Z}Ge_{\mathbb{Q}}(\psi^{G})))$. Furthermore, it can be  checked that $\mathfrak{C}$ is a subgroup. We next see that the index of $\mathfrak{C}$ in $\mathcal{Z}(\mathcal{U}(\mathbb{Z}(1-e_{\mathbb{Q}}(\psi^{G}))+\mathbb{Z}Ge_{\mathbb{Q}}(\psi^{G})))$ is finite. Consider an arbitrary $w \in \mathcal{Z}(\mathcal{U}(\mathbb{Z}(1-e_{\mathbb{Q}}(\psi^{G}))+\mathbb{Z}Ge_{\mathbb{Q}}(\psi^{G})))$. We can write $w = \pm (1-e_{\mathbb{Q}}(\psi^{G}))+w'$, where $w' \in \mathcal{Z}(\mathcal{U}(\mathbb{Z}Ge_{\mathbb{Q}}(\psi^{G})))$. W.l.o.g. we can assume that $w = 1-e_{\mathbb{Q}}(\psi^{G})+w'$. By Proposition \ref{p2}, $w'=\sum_{t \in T}w''^{t}$, where $w''$ is a central unit of $\mathbb{Z}Ce_{\mathbb{Q}}(\psi)$. As $\psi$ has the least center property w.r.t. $G$, $\mathcal{Z}(\mathbb{Q}Ce_{\mathbb{Q}}(\psi))= \mathcal{Z}(\mathbb{Q}Se_{\mathbb{Q}}(\psi))^{\operatorname{Cen}_{G}(e_{\mathbb{Q}}(\psi))/S}$. Hence $w'' \in \mathcal{Z}(\mathbb{Z}Se_{\mathbb{Q}}(\psi))^{\operatorname{Cen}_{G}(e_{\mathbb{Q}}(\psi))/S}$. Now, $w = 1-e_{\mathbb{Q}}(\psi^{G})+w' = 1-e_{\mathbb{Q}}(\psi^{G})+ \sum_{t \in T}w''^{t} = \prod_{t \in T} (1-e_{\mathbb{Q}}(\psi)+w'')^{t}$. Since the entries in the product commute and a power of $1-e_{\mathbb{Q}}(\psi)+w''$ belongs to $\mathfrak{B}$ (by step I), it follows that a power of $w$ belongs to $\mathfrak{C}$. \par It now remains to show that the elements of $\mathfrak{C}$ are central units of $\mathbb{Z}G$. Observe that an arbitrary element of $\mathcal{Z}(\mathcal{U}(\mathbb{Z}(1-e_{\mathbb{Q}}(\psi^{G}))+\mathbb{Z}Ge_{\mathbb{Q}}(\psi^{G})))$ commutes with all $g \in G$. Therefore, $\mathfrak{C}$ is contained in the center of $\mathbb{Z}G$. Furthermore, the elements of $A$ being units in $\mathbb{Z}S$, it follows that the elements of $\mathfrak{C}$ are units in $\mathbb{Z}G$. This proves step II and completes the proof of the lemma. $~\Box$\vspace{.25cm} \para Let $(H,K)$ be a generalized strong Shoda pair of $G$ and let $\mathcal{N}:~H=H_{0}\leq H_{1}\leq \cdots \leq H_{n}=G$ be a strong inductive chain from $H$ to $G$. Let $\lambda \in \operatorname{Lin}(H,K)$. If $\lambda^{H_{i}}$ has the least center property w.r.t. $H_{i+1}$ for all $0\leq i \leq n-1$, then we say that the \textit{generalized strong Shoda pair} $(H,K)$ \textit{of} $G$ \textit{has the least center property w.r.t.} $\mathcal{N}$. Note that this definition is independent of the choice of $\lambda$.\par Consider a generalized strong Shoda pair $(H,K)$ of $G$ which  has the least center property  w.r.t. the strong inductive chain $\mathcal{N}:~H=H_{0}\leq H_{1}\leq \cdots \leq H_{n}=G$. We'll use Lemma \ref{l5} to provide an iterative process to construct central units of $\mathbb{Z}G$ from the central units of $\mathbb{Z}H$ lying in $\mathcal{Z}(\mathcal{U}(\mathbb{Z}(1-\varepsilon(H,K))+\mathbb{Z}H\varepsilon(H,K)))$. Let $\lambda \in \operatorname{Lin}(H,K)$. For $u \in \mathcal{Z}(\mathcal{U}(\mathbb{Z}H))\cap \mathcal{Z}(\mathcal{U}(\mathbb{Z}(1-\varepsilon(H,K))+\mathbb{Z}H\varepsilon(H,K)))$, put $$z_{0}^{\mathcal{N}}(u)=u$$ and for $0\leq i\leq n-1$, put $$z_{i+1}^{\mathcal{N}}(u)=\prod_{t \in T_{i}}\big(\prod_{c \in C_{i}}z_{i}^{\mathcal{N}}(u)^{c}\big)^{t},$$ where $C_{i}=\operatorname{Cen}_{H_{i+1}}(e_{\mathbb{Q}}(\lambda^{H_{i}}))$ and $T_{i}$ is a right transversal of $C_{i}$ in $H_{i+1}$. Denote the final step of the construction $z_{n}^{\mathcal{N}}(u)$ by $z^{\mathcal{N}}(u)$. From Lemma \ref{l5}, it follows that the construction of $z_{i}^{\mathcal{N}}(u)$ is well defined and it belongs to $\mathcal{Z}(\mathcal{U}(\mathbb{Z}H_{i}))\cap\mathcal{Z}(\mathcal{U}(\mathbb{Z}(1-e_{\mathbb{Q}}(\lambda^{H_{i}}))+\mathbb{Z}H_{i}e_{\mathbb{Q}}(\lambda^{H_{i}})))$ for all $0\leq i\leq n$. In particular, $z^{\mathcal{N}}(u)$ is a central unit of $\mathbb{Z}G$.\par Given $g \in G$ and $k,m$ positive integers such that $k^{m}\equiv 1\operatorname{mod}|g|$, where $|g|$ is the order of $g$, the following is a unit of $\mathbb{Z}G$: $$u_{k,m}(g)=(1+g+\cdots +g^{k-1})^{m}+\frac{1-k^{m}}{|g|}(1+g+\cdots +g^{|g|-1}).$$ The units of this form are called Bass cyclic units based on $g$ with parameters $k$ and $m$ and were introduced by Bass \cite{B}. \par When $M$ is a normal subgroup of $G$, then $$u_{k,m}(1-\widehat{M}+g\widehat{M})=1-\widehat{M}+u_{k,m}(g)\widehat{M}$$ is an invertible element of $\mathbb{Z}G(1-\widehat{M})+\mathbb{Z}G\widehat{M}$. As this is an order in $\mathbb{Q}G$, for each element $b=u_{k,m}(1-\widehat{M}+g\widehat{M})$ there is a positive integer $n$ such that $b^{n} \in \mathcal{U}(\mathbb{Z}G)$. Let $n_{b}$ denote the minimal positive integer satisfying this condition. The element $$u_{k,m}(1-\widehat{M}+g\widehat{M})^{n_{b}}=1-\widehat{M}+u_{k,mn_{b}}(g)\widehat{M}$$ is called the generalized Bass unit based on $g$ and $M$ with parameters $k$ and $m$. These units were introduced by Jespers and Parmenter in \cite{JP}. Note that if $M$ is such that $G' \subset M$, then the generalized Bass units are central units in $\mathbb{Z}G$.\begin{theorem}\label{t4} Let $G$ be a generalized strongly monomial group and let \linebreak $\mathcal{S}=\{(H_{i},K_{i})~|~1\leq i \leq n\}$ be a complete and irredundant set of generalized strong Shoda pairs of $G$. For each $i$, let $A_{(H_{i},K_{i})}$ be the subgroup of $\mathcal{Z}(\mathcal{U}(\mathbb{Z}H_{i}))$ generated by the generalized Bass units $b_{i}^{n_{b_{i}}}$, where $b_{i}=u_{k,m}(1-\widehat{H_{i}'}+h\widehat{H_{i}'})$, $h \in H_{i}$, $k$ and $m$ positive integers s.t. $k^{m}\equiv 1\operatorname{mod}|h|$. If $(H_{i},K_{i})$ has the least center property w.r.t. a strong inductive chain $\mathcal{N}_{i}$ from $H_{i}$ to $G$ for all $1\leq i \leq n$, then $$\small\{z^{\mathcal{N}_{1}}(u_{1})\cdots z^{\mathcal{N}_{n}}(u_{n})~|~u_{i}\in A_{(H_{i},K_{i})}\cap \mathcal{Z}(\mathcal{U}(\mathbb{Z}(1-\varepsilon(H_{i},K_{i}))+\mathbb{Z}H_{i}\varepsilon(H_{i},K_{i}))),~1\leq i\leq n\}$$ forms a subgroup of $\mathcal{Z}(\mathcal{U}(\mathbb{Z}G))$ which is contained in the group generated by  generalized Bass units of $\mathbb{Z}G$ and its index in $\mathcal{Z}(\mathcal{U}(\mathbb{Z}G))$ is finite.\end{theorem}{\bf Proof.} For $1\leq i\leq n$, fix a $\lambda_{i} \in \operatorname{Lin}(H_{i},K_{i})$. Using Bass Milnor theorem, it follows that $A_{(H_{i},K_{i})}$ is a subgroup of finite index in $\mathcal{U}(\mathbb{Z}(1-\widehat{H_{i}'})+\mathbb{Z}H_{i}\widehat{H_{i}'})$. Let $\mathcal{O}_{1}^{(i)}=\mathbb{Z}(1-\widehat{H_{i}'})+\mathbb{Z}H_{i}\widehat{H_{i}'},~\mathcal{O}_{2}^{(i)}=\mathbb{Z}(1-\varepsilon(H_{i},K_{i}))+\mathbb{Z}H_{i}
  \varepsilon(H_{i},K_{i}),~B_{1}^{(i)}=\mathbb{Q}(1-\widehat{H_{i}'})+\mathbb{Q}H_{i}\widehat{H_{i}'}$ and $B_{2}^{(i)}=\mathbb{Q}(1-\varepsilon(H_{i},K_{i}))+\mathbb{Q}
  H_{i}\varepsilon(H_{i},K_{i}).$ Clearly $\mathcal{O}_{1}^{(i)}$ and $\mathcal{O}_{2}^{(i)}$ are $\mathbb{Z}$-orders in $B_{1}^{(i)}$ and $B_{2}^{(i)},$ respectively. Moreover $B_{2}^{(i)}\subseteq B_{1}^{(i)}.$ So using Lemma 4.6.6 of \cite{JdR}, it immediately follows that $\mathcal{O}_{1}^{(i)}\cap B_{2}^{(i)}$ is a $\mathbb{Z}$-order in $B_{2}^{(i)}$. This, in view of (\cite{JdR}, Lemma 4.6.9) gives that $\mathcal{U}(\mathcal{O}_{1}^{(i)})\cap \mathcal{U}(\mathcal{O}_{2}^{(i)})$ is of finite index in $\mathcal{U}(\mathcal{O}_{2}^{(i)})$. Consequently $A_{(H_{i},K_{i})}\cap \mathcal{U}(\mathcal{O}_{2}^{(i)})$ is of finite index in $\mathcal{U}(\mathcal{O}_{2}^{(i)})$. Observe that $\mathcal{O}_{2}^{(i)}$ is commutative and hence $\mathcal{U}(\mathcal{O}_{2}^{(i)}) = \mathcal{Z}(\mathcal{U}(\mathcal{O}_{2}^{(i)}))$. We have thus obtained the subgroup \linebreak $A_{(H_{i},K_{i})}\cap \mathcal{Z}(\mathcal{U}(\mathcal{O}_{2}^{(i)}))$ of $\mathcal{Z}(\mathcal{U}(\mathbb{Z}H_{i}))\cap \mathcal{Z}(\mathcal{U}(\mathbb{Z}(1-\varepsilon(H_{i},K_{i}))+\mathbb{Z}H_{i}\varepsilon(H_{i},K_{i})))$ such that its index in $\mathcal{Z}(\mathcal{U}(\mathbb{Z}(1-\varepsilon(H_{i},K_{i}))+\mathbb{Z}H_{i}\varepsilon(H_{i},K_{i})))$ is finite. By the repeated application of Lemma \ref{l5}, we thus obtain that $$U_{i}=\{ z^{\mathcal{N}_{i}}(u_{i})~|~u_{i} \in A_{(H_{i},K_{i})}\cap\mathcal{Z}(\mathcal{U}(\mathcal{O}_{2}^{(i)}))\}$$ is a subgroup of $\mathcal{Z}(\mathcal{U}(\mathbb{Z}G))$ and its index in $\mathcal{Z}(\mathcal{U}(\mathbb{Z}(1-e_{\mathbb{Q}}(\lambda_{i}^{G}))+\mathbb{Z}Ge_{\mathbb{Q}}(\lambda_{i}^{G})))$ is finite. Now, consider the following subgroup $$U=\{z^{\mathcal{N}_{1}}(u_{1})z^{\mathcal{N}_{1}}(u_{2})\cdots z^{\mathcal{N}_{n}}(u_{n})~|~u_{i}\in A_{(H_{i},K_{i})}\cap \mathcal{Z}(\mathcal{U}(\mathcal{O}_{2}^{(i)})),~1\leq i\leq n\}$$ of $\mathcal{Z}(\mathcal{U}(\mathbb{Z}G))$. We will show that the index of $U$ in $\mathcal{Z}(\mathcal{U}(\mathbb{Z}G))$ is finite. For this, it is enough to see that a power of each element of $\mathcal{Z}(\mathcal{U}(\mathbb{Z}G))$ belongs to $U$. Consider an arbitrary $u \in \mathcal{Z}(\mathcal{U}(\mathbb{Z}G))$. We can write this element as $$u=\sum_{1\leq i\leq n}ue_{\mathbb{Q}}(\lambda_{i}^{G})=\prod_{1\leq i\leq n}(1-e_{\mathbb{Q}}(\lambda_{i}^{G})+ue_{\mathbb{Q}}(\lambda_{i}^{G})).$$ As, for any $i$,  a power of $1-e_{\mathbb{Q}}(\lambda_{i}^{G})+ue_{\mathbb{Q}}(\lambda_{i}^{G})$ belongs to $U_{i}$ and the entries in the product on the right side commute, it follows that a power of $u$ belongs to $U$.\par Finally, we can see that $U$ is contained in the group generated by the generalized Bass units because the conjugates of Bass units are again Bass units and $(H^{g},K^{g})$ is a generalized strong Shoda pair of $G$ if $(H,K)$ is a generalized strong Shoda pair of $G$. This completes the proof. $~\Box$\vspace{.2cm} \para In the above theorem, the assumption that $(H_{i},K_{i})$ has the least center property w.r.t. $\mathcal{N}_{i}$ is satisfied when $(H_{i},K_{i})$ is a strong Shoda pair of $G$ and $\mathcal{N}_{i}$ is the strong inductive chain $H_{i} \leq G$. Thus the above theorem is applicable to strongly monomial groups and we have the following: \begin{cor}\label{c1} {\rm (\cite{JOdRVG}, Theorem 5.1)} If $G$ is a strongly monomial group, then the group generated by generalized Bass units contain a subgroup of finite index in $\mathcal{Z}(\mathcal{U}(\mathbb{Z}G))$.\end{cor}\vspace{.25cm}\par We now show that the subgroup of finite index in $\mathcal{Z}(\mathcal{U}(\mathbb{Z}G))$ given in the  Theorem \ref{t4} can be simplified when $G$ is a subnormally monomial group. If $g\in G$ is such that $\langle g \rangle$ is subnormal in $G$, then a construction of central unit of $\mathbb{Z}G$ from the Bass unit based on $g$ is given in \cite{JOdRVG}. The same idea can be used to provide a construction of central units of $\mathbb{Z}G$ from central units of $\mathbb{Z}H$, when $H$ is a subnormal subgroup of $G$. Let $\mathcal{N}:H=H_{0}\unlhd H_{1}\unlhd \cdots \unlhd H_{n}=G$ be a subnormal series. For $u \in \mathcal{Z}(\mathcal{U}(\mathbb{Z}H))$, put $c_{0}^{\mathcal{N}}(u)=u$ and for $0\leq i\leq n-1$, set $$c_{i+1}^{\mathcal{N}}(u)=\prod_{t \in T_{i}}c_{i}^{\mathcal{N}}(u)^{t},$$ where $T_{i}$ is a right transversal of $H_{i}$ in $H_{i+1}$. Denote $c_{n}^{\mathcal{N}}(u)$ by $c^{\mathcal{N}}(u)$. Analogous to (\cite{JOdRVG}, Lemma 3.1), the following properties hold:\begin{description}\item [(i)] $c_{i}^{\mathcal{N}}(u)^{x}\in \mathbb{Z}H_{i}$ for all $x \in H_{i+1}$;\item [(ii)] $c_{i}^{\mathcal{N}}(u)^{x}=c_{i}^{\mathcal{N}}(u)$ for all $x \in H_{i}$ ;\item [(iii)]  $c_{i}^{\mathcal{N}}(u)$ is independent of the transversal $T_{i-1}$.\end{description} In particular, we conclude that $c^{\mathcal{N}}(u)$ is well defined and is a central unit of $\mathbb{Z}G$. \begin{cor}\label{c2} Let $G$ be a subnormally monomial group. Let $\{(H_{i},K_{i})~|~1\leq i\leq n\}$ be a complete and irredundant set of Shoda pairs of $G$, where each $H_{i}$ is subnormal in $G$. For each $i$, fix a subnormal series $\mathcal{N}_{i}$ from $H_{i}$ to $G$. Let $A_{(H_{i},K_{i})}$ be as in Theorem \ref{t4}. If $(H_{i},K_{i})$ has the least center property w.r.t. $\mathcal{N}_{i}$ for all $1\leq i\leq n$, then $$\{c^{\mathcal{N}_{1}}(u_{1})\cdots c^{\mathcal{N}_{n}}(u_{n})~|~u_{i}\in A_{(H_{i},K_{i})}\cap \mathcal{Z}(\mathcal{U}(\mathbb{Z}(1-\varepsilon(H_{i},K_{i}))+\mathbb{Z}H_{i}\varepsilon(H_{i},K_{i}))),~1\leq i\leq n\}$$ is a subgroup of finite index in $\mathcal{Z}(\mathcal{U}(\mathbb{Z}G))$.\end{cor}{\bf Proof.} Observe that any subnormal series from $H_{i}$ to $G$ is a strong inductive chain from $H_{i}$ to $G$. In particular, $\mathcal{N}_{i}$ is a  strong inductive chain from $H_{i}$ to $G$. Also note that a power of  $c^{\mathcal{N}_{i}}(u_{i})$ is $z^{\mathcal{N}_{i}}(u_{i})$ for all $ 1 \leq i \leq n$. Therefore, Theorem \ref{t4} yields the result. $\Box$ \vspace{.25cm}\para If $G$ is an abelian-by-supersolvable group with the property that every cyclic subgroup of order not a divisor of 4 or 6 is subnormal in $G$, then it is known (\cite{JOdRVG}, Theorem 3.2) that the group generated by Bass units of $\mathbb{Z}G$ contains a subgroup of finite index in $\mathcal{Z}(\mathcal{U}(\mathbb{Z}G)).$  Generalizing this result, we prove the following:\begin{theorem}\label{t5} Let $G \in \mathcal{C}$ be such that every cyclic subgroup of order not a divisor of $4$ or $6$ is subnormal in $G$. Let $S$ be a complete and irredundant set of generalized strong Shoda pairs of $G$. If each $(H,K) \in S$ has the least center property w.r.t. some strong inductive chain from $H$ to $G$, then, the group generated by Bass units of $\mathbb{Z}G$ contains a subgroup of finite index in $\mathcal{Z}(\mathcal{U}(\mathbb{Z}G))$.\end{theorem} \par We need the following lemma. We'll just state it without proof as it can be obtained by working as in the first part of the proof of (\cite{JOdRVG}, Theorem 3.2). \begin{lemma}\label{l6} Let $G$ be a finite group such that each cyclic subgroup of order not a divisor of 4 or 6 is subnormal in $G$. For each such cyclic subgroup $\langle g \rangle$, fix a subnormal series $\mathcal{N}_{g}$ from $\langle g \rangle$ to $G$. Then, given a unit $z$ in $\mathbb{Z}G\widehat{G'}$ there exist Bass units $b_{g_{1}},b_{g_{2}},\cdots,b_{g_{r}}$ in $\mathbb{Z}G$ such that $\prod_{1\leq i\leq r}c^{\mathcal{N}_{g_{i}}}(b_{g_{i}})\widehat{G'}$ is a power of $z$.\end{lemma} {\bf Proof of Theorem \ref{t5}.} The theorem is proved once we show that a power of each central unit of $\mathbb{Z}G$ is a product of Bass units of $\mathbb{Z}G$, as $\mathcal{Z}(\mathcal{U}(\mathbb{Z}G))$ is a finitely generated abelian group. We'll prove it by induction  on $|G|$. If $|G|=1$, the result is trivially true. Assume the result holds for the groups of order strictly less than $|G|$. We can assume that $G$ is non abelian, as for abelian groups the result follows from Bass Milnor Theorem. Consider an arbitrary $z \in \mathcal{Z}(\mathcal{U}(\mathbb{Z}G))$. We have to show that a power of $z$ is a product of Bass units of $\mathbb{Z}G$. We can write $z=z'+z''$, where $z' \in \mathcal{Z}(\mathcal{U}(\mathbb{Z}G(1-\widehat{G'})))$ and $z'' \in \mathcal{U}(\mathbb{Z}G\widehat{G'})$. By Lemma \ref{l6}, there exist Bass units $b_{g_{1}},b_{g_{2}},\cdots,b_{g_{r}}$ in $\mathbb{Z}G$ such that $c^{\mathcal{N}_{g_{1}}}(b_{g_{1}})\cdots c^{\mathcal{N}_{g_{r}}}(b_{g_{r}})\widehat{G'}=z''^{m}$ for some integer $m\geq 1$. Hence $$z^{m}(\prod_{1\leq i\leq r}c^{\mathcal{N}_{g_{i}}}(b_{g_{i}}))^{-1}=z'''+\widehat{G'}~{\rm for~some}~ z''' \in \mathcal{Z}(\mathcal{U}(\mathbb{Z}G(1-\widehat{G'}))).$$ Observe that the central unit $c^{\mathcal{N}_{g_{i}}}(b_{g_{i}})$ of $\mathbb{Z}G$ is a product of Bass units of $\mathbb{Z}G$ for all $1\leq i\leq r.$ Therefore, it is enough to show that a power of $z'''+\widehat{G'}$ is a product of Bass units of $\mathbb{Z}G$. For $(H,K)\in S$, fix a linear character $\lambda_{(H,K)}\in \operatorname{Lin}(H,K)$. We have $$1-\widehat{G'}=\sum_{\stackrel{(H,K)\in S}{H\neq G}}e_{\mathbb{Q}}(\lambda_{(H,K)}^{G}).$$ Hence $$\displaystyle z'''+\widehat{G'}=\displaystyle\sum_{\stackrel{(H,K)\in S}{H\neq G}}z'''e_{\mathbb{Q}}(\lambda_{(H,K)}^{G})+\widehat{G'}=\displaystyle\prod_{\stackrel{(H,K)\in S}{H\neq G}} (1-e_{\mathbb{Q}}(\lambda_{(H,K)}^{G})+z'''e_{\mathbb{Q}}(\lambda_{(H,K)}^{G})).$$ As the entries in the product on the right side commute, it suffices to show that a power of $1-e_{\mathbb{Q}}(\lambda_{(H,K)}^{G})+z'''e_{\mathbb{Q}}(\lambda_{(H,K)}^{G})$ is a product of Bass units of $\mathbb{Z}G$ for every $(H,K) \in S$ with $H\neq G$. Consider an arbitrary $(H,K) \in S$ with $H \neq G$. As the class $\mathcal{C}$ is closed under subgroups, we have $H \in \mathcal{C}$. Since $|H|<|G|,$ by induction there exists a subgroup $A$ in the group generated by Bass units of $\mathbb{Z}H$ such that $A$ is of finite index in $\mathcal{Z}(\mathcal{U}(\mathbb{Z}H))$. As a consequence, it can be shown that $A(H,K)=A\cap (\mathbb{Z}(1-\varepsilon(H,K))+\mathbb{Z}H\varepsilon(H,K))$ is of finite index in $\mathcal{U}(\mathbb{Z}(1-\varepsilon(H,K))+\mathbb{Z}H\varepsilon(H,K))$. This yields the subgroup $A(H,K)$ of $\mathcal{Z}(\mathcal{U}(\mathbb{Z}H))\cap
 \mathcal{Z}(\mathcal{U}(\mathbb{Z}(1-\varepsilon(H,K))+\mathbb{Z}H\varepsilon(H,K)))$ such that its index in $\mathcal{Z}(\mathcal{U}(\mathbb{Z}(1-\varepsilon(H,K))+\mathbb{Z}H\varepsilon(H,K)))$ is finite. Thus applying Lemma \ref{l5} repeatedly, we obtain that $B(H,K)=\{z^{\mathcal{N}}(u)~|~u \in A(H,K)\}$ is a subgroup of $\mathcal{Z}(\mathcal{U}(\mathbb{Z}G))$ and its index in $\mathcal{Z}(\mathcal{U}(\mathbb{Z}(1-e_{\mathbb{Q}}(\lambda_{(H,K)}^{G}))+\mathbb{Z}Ge_{\mathbb{Q}}(\lambda_{(H,K)}^{G})))$ is finite, where $\mathcal{N}$ is a strong inductive chain w.r.t. $G$ such that $(H,K)$ has the least center property w.r.t. $\mathcal{N}$. This proves the desired result as the elements of $B(H,K)$ are product of Bass units of $\mathbb{Z}G$. $~\Box$\vspace{.25cm}\\ \textbf{Remark:} In Theorem \ref{t5}, one can further show that $$ \langle c^{\mathcal{N}_{g}}(b_{g})~|~b_{g}~{\rm is~a~Bass~unit~based~on~}g,~g \in G \rangle $$ is a subgroup of the group generated by Bass units of $\mathbb{Z}G$ which is of finite index in $\mathcal{Z}(\mathcal{U}(\mathbb{Z}G)),$ where $\mathcal{N}_{g}$ is a subnormal series from the cyclic subgroup $\langle g \rangle$ to $G$. This can be seen by proceeding as in (\cite{JOdRVG}, Corollary 3.4). Moreover, one can construct a virtual basis of $\mathcal{Z}(\mathcal{U}(\mathbb{Z}G))$ consisting of Bass units of $\mathbb{Z}G$ using the technique followed in (\cite{JOdRVG}, Theorem 4.2). \bibliographystyle{amsplain}
\bibliography{BIB1}

\end{document}